\documentclass[12pt,letterpaper]{amsart}
\usepackage{amssymb,amsmath}

\newtheorem{thm}{Theorem}
\newtheorem{prop}[thm]{Proposition}

\newcommand{\del}[1]{\frac{\partial}{\partial #1}}

\begin{document}

\title[The geometry of Chazy's equations]{The geometry of Chazy's homogeneous third-order differential equations}
\author{Adolfo GUILLOT}
\address{Instituto de Matem\'aticas, Unidad Cuernavaca,
 Universidad Nacional Aut\'onoma de M\'exico,
A.P.~273-3 Admon.~3,
Cuernavaca, Morelos, 62251, Mexico.}
\email{adolfo@matcuer.unam.mx}

\thanks{The author was partially supported by CONACyT grant 58354 and PAPIIT-UNAM grant IN104510.}%
\thanks{MSC 2000 34M55, 34M45, 32M25}%

\begin{abstract}
Chazy studied a family of homogeneous third-order autonomous differential equations. They are those, within a certain class, admitting exclusively single-valued solutions. Each one of these equations yields a polynomial vector field in  complex three-dimensional space. For almost all of these these vector fields, the Zariski closure of a generic orbit yields an affine surface endowed with  a holomorphic vector field that has exclusively single-valued solutions. We classify these surfaces and relate this classification to recent results of Rebelo and the author.
\end{abstract}
\maketitle

\section{Introduction}
In his attempt to extend Painlev\'e's theory of second order differential equations in the complex domain to  equations of the third order, Chazy, in~\cite{chazy}, studied those  differential equations that have only single-valued solutions and that are of the form
$\phi'''=P(\phi'',\phi',\phi)$,
where~$P$ is a polynomial whose coefficients are meromorphic functions in the independent variable. In dealing with such a problem, a first step is to determine the differential equations having only single-valued solutions that have the same form, with~$P$ a polynomial with \emph{constant} coefficients and some form of homogeneity. An instance of this simplified problem is given by the equations of the form
\begin{equation}\label{eq:chazy}\phi'''=a_3\phi^4+a_2\phi^2\phi'+a_1(\phi')^2+\delta\phi\phi''\;\;\;\; (a_1,a_2,a_3,\delta\in\mathbf{C}).\end{equation}
Chazy found~\cite[p.~335--337]{chazy} that the equations of the form~(\ref{eq:chazy}) having only single-valued solutions belong, up to a natural rescaling, to the following list:
\begin{align}
\tag{II} \phi''' & = 2\phi\phi''+2(\phi')^2 \displaybreak[0]\\
\tag{IV} \phi''' & =  3\phi\phi''+3(\phi')^2 -3\phi^2\phi' \displaybreak[0]\\
\tag{V} \phi''' & = 2\phi\phi''+4(\phi')^2-2\phi^2\phi'  \displaybreak[0]\\
\tag{VI} \phi''' & = \phi\phi''+5(\phi')^2-\phi^2\phi' \displaybreak[0]\\
\tag{VII} \phi''' & = \phi\phi''+ 2\phi^2\phi'+2(\phi')^2 \displaybreak[0]\\
\tag{VIII} \phi''' & =  6\phi^2\phi' \displaybreak[0]\\
\tag{IX} \phi''' & = 18(\phi'+\phi^2)(\phi'+3\phi^2)-6(\phi')^2 \displaybreak[0]\\
\tag{X} \phi''' & =  6\phi^2\phi'+3\frac{9+7\sqrt{3}}{11}(\phi'+\phi^2)^2 \displaybreak[0]\\
\tag{XI} \phi''' & = {\textstyle\frac{1}{2}}\kappa \phi\phi''+\left({\textstyle\frac{1}{2}}\kappa+6\right)(\phi')^2-3\kappa\phi^2\phi'+{\textstyle\frac{3}{8}}\kappa^2\phi^4,\;\kappa=1-k^2, k\in\mathbf{Z}\setminus 6\mathbf{Z} \\
\tag{XII} \phi''' & =  2 \phi\phi''-3(\phi')^2+\frac{4}{36-k^2}(6\phi'-\phi^2)^2,\;k\in \mathbf{Z}\cup\{\infty\},\;  k\neq 0,1,6;
\end{align}
(there are two inequivalent equations behind Equation~X, one for each  determination of~$\sqrt{3}$). By setting~$(x,y,z)=(\phi,\phi',\phi'')$, equation~(\ref{eq:chazy}) may be written as the  vector field
\begin{equation}\label{field:chazy}V=y\frac{\partial}{\partial x}+z\frac{\partial}{\partial y}+(a_3x^4+a_2x^2y+ a_1y^2+\delta xz)\frac{\partial}{\partial z}\end{equation}
in~$\mathbf{C}^3$. The Zariski closure of an orbit is an affine variety~$\Sigma\subset\mathbf{C}^3$, endowed with the restriction of~$V$ to it. If the original differential equation has only single-valued solutions, this vector field will have only single valued solutions, this is, it will be \emph{semicomplete} in the sense of Rebelo~\cite{rebelo}. When the dimension of this variety is strictly smaller than three, we have a fairly good understanding of the situation. In dimension one, the variety embeds into a rational or elliptic curve, where the vector field extends as a holomorphic one. In dimension two, we have the following result from~\cite{gui-reb}:

\begin{thm} \label{guillot-rebelo} Let~$\Sigma$ be a (possibly singular) affine algebraic surface in~$\mathbf{C}^n$ and~$V$ the restriction of a polynomial vector field in~$\mathbf{C}^n$ to~$\Sigma$. If~$V$ is semicomplete in restriction to~$\Sigma$ then, up to a birational transformation, there is a compact projective surface~$S$ endowed with a meromorphic vector field~$W$ and an embedding~$j:\Sigma\to S$, with Zariski open image, mapping~$V$ to~$W$, such that either
\begin{itemize}\item $W$ is a holomorphic vector field or
 \item $S$ is a rational or elliptic fibration (over a rational or elliptic base), the poles of~$W$ are a union of fibers and, at the other fibers, $W$~is either transverse or tangent.
\end{itemize}
 \end{thm}

Holomorphic vector fields in complex projective surfaces, related (after desingularization) to the first item of this Theorem, are very well understood~\cite{kob}. A holomorphic vector field on a complex projective surface may or may not have a meromorphic first integral. If it does, the Zariski closure of an orbit is a curve: the surface is an elliptic or rational fibration and the vector field tangent to the fibration. In the absence of a first integral, we have the following result (see, for example, \cite{bru}):
\begin{prop} A holomorphic vector field in a compact complex projective surface that does not have a meromorphic first integral is either:
\begin{itemize}
 \item a vector field on an Abelian surface,
 \item a vector field without zeroes on a rational fibration over an elliptic curve or
\item a holomorphic vector field on a rational surface.
\end{itemize}
\end{prop}

As we shall see, besides the equations Chazy~XII for~$k\geq 7$, all the equations in Chazy's list have at least a polynomial first integral and thus, in restriction to a level surface, we are exactly in the setting of Theorem~\ref{guillot-rebelo}. The aim of this article is to describe explicitly the geometry behind the Chazy equations, by identifying in which of the possibilities given by Theorem~\ref{guillot-rebelo} does the corresponding equation fit. Our results can be summarized as follows:

\begin{thm} For the vector fields corresponding to Chazy's equations (except Chazy~XII for~$k\geq 7$), the Zariski closure of the generic orbit is
\begin{itemize}
\item a curve, for
\begin{itemize}
\item Chazy~XII for~$k\in\{2,3,4,5\}$: the generic orbit embeds into  a rational curve.
\item Chazy~VII and Chazy~VIII : the generic orbit embeds into an elliptic curve.
\end{itemize}
\item a surface that can be compactified as
\begin{itemize}
\item a rational fibration over an elliptic curve for Chazy~XI, where the vector field extends as a holomorphic one (transverse to the fibers);
\item a rational fibration over a rational curve in Chazy~II;
 \item a rational surface for Chazy~IV, Chazy~V and Chazy~VI, where the vector field extends as a holomorphic one;
\item an Abelian surface for Chazy~IX and~Chazy~X, where the vector field extends as a (nowhere-zero) holomorphic one.
\end{itemize}
\end{itemize}
\end{thm}

Our results will show explicitly how the above compactifications are, and this involves, at times, some very complicated formulas. Our calculations were made with the help of a computer algebra system, but the information we include is sufficient to rebuild our calculations. We arrived at our results by systematically looking for quasihomogeneous first integrals, quasihomogeneous  invariant hypersurfaces and quasihomogenous   vector fields commuting with the corresponding vector field~(\ref{field:chazy}), except in equations~IX and~X where we took Cosgrove's solutions of these equations~\cite{cosgrove-2} as a starting point. In Section~\ref{cuentas} we detail, in a more or less self-contained fashion, our approach to the classical procedure that, starting from the equations of the form~(\ref{eq:chazy}), leads to Chazy's list. \\

Chazy's own list contains, beyond the previous ones, the equations
\begin{align}
\tag{I} \phi''' & = -6(\phi')^2, \\
\tag{III} \phi''' & =  2 \phi\phi''-3(\phi')^2,\\
\tag{XIII} \phi''' & =  12\phi\phi'.
\end{align}
Equation~I is equation~XI when~$k=1$ (up to an inessential change of sign) and equation~III is equation~XII when~$k=\infty$. We have not hesitated to include these equations within the larger families. Equation~XIII does not belong to the family~(\ref{eq:chazy}).\\

In this article we use the property of \emph{semicompleteness} to formalize the notion of ``complex differential equations without multivalued solutions''. There are other more or less related notions that allow us to talk of the presence and nature of multivaluedness phenomena within the solutions of complex differential equations, like the \emph{Painlev\'e Property} or the \emph{absence of movable critical points}. Semicompleteness is a property that is either present or absent in meromorphic vector fields on complex manifolds \cite{gui-reb} and was firstly defined for  holomorphic vector fields~\cite{rebelo}. It is particularly well-adapted to the geometric setting that will concern us. Its presence/absence is a birational invariant, complete vector fields (like vector fields on compact manifolds) are semicomplete and, furthermore, the restriction of a semicomplete vector field on a manifold to an open subset will remain semicomplete (these are the properties that will be actually used in this article).

\section{The geometry of Chazy's equations}
We will go through all the equations (II)-(XII) in order to describe their geometry.
We will denote by~$V$ the corresponding vector field in~$\mathbf{C}^3$ of the form~(\ref{field:chazy}). The fact that these vector fields are semicomplete will follow from our geometric description. The key point is that, with respect to the linear vector field
\begin{equation}\label{ell}L=x\frac{\partial}{\partial x}+2y\frac{\partial}{\partial y}+3z\frac{\partial}{\partial z},\end{equation}
the vector field~$V$ satisfies the Lie bracket relation~$[L,V]=V$. In other words, if we give the variables $x$, $y$ and~$z$ the weights~$1$, $2$ and~$3$, the vector fields become homogeneous. As we shall see, except for Chazy~XII with~$k\geq 7$, all the vector fields admit quasihomogeneous polynomial first integrals, this is, polynomials~$Q$ such that~$L\cdot Q=(d-1)Q$ for some integer~$d$ (the degree of the quasihomogeneous polynomial) and such that~$V\cdot Q=0$  ($V\cdot Q$ stands for the derivative of~$Q$ along~$V$). We will analyze the vector field in restriction to a single non-zero level surface of this first integral (by the quasihomogeneity of the vector field and of the first integral, all the nonzero level surfaces are equivalent).
As a consequence of the quasihomogeneity of the first integral, any level surface will be preserved by the order~$d$ cyclic group of transformations generated by
$$(x,y,z)\mapsto (\rho x,\rho^2y,\rho^3z)$$
for~$\rho$ a primitive~$d^\text{th}$ root of unity, and under this transformation, the vector field will be multiplied by a constant. We will also describe these symmetries.

We will not study the vector field in restriction to the surface where the first integral vanishes (semicompleteness in restriction to this surface follows from the semicompleteness of the vector field in its complement).

\subsection*{Chazy II} \emph{There is a quasihomogeneous polynomial first integral. The generic level surface compactifies as a rational fibration over a rational curve and the vector field extends as a meromorphic one.}
In the quasihomogeneous coordinates
$(X,Y,Z)=(x,y-x^2,z-2xy)$,
the vector field~$V$ becomes
$$(X^2+Y)\del{X}+Z\del{Y}.$$
It has the first integral~$Z$. We can embed~$\Sigma=Z^{-1}(1)$ into~$\mathbf{CP}^1\times \mathbf{CP}^1$ by~$(X,Y)\mapsto([X:1],[Y:1])$. The vector field~$V|_\Sigma$ extends to a meromorphic vector field in~$\mathbf{CP}^1\times \mathbf{CP}^1$. Under the projection to the second factor we get the vector field~$\partial/\partial Y$. Any solution to this vector field may be lifted to the whole space and the vector field is thus semicomplete (in terms of functions, $X(t)$ is a solution to the Riccati equation~$\phi'=\phi^2+t$).

\subsection*{Chazy IV} \emph{There is a quasihomogeneous polynomial first integral. The restriction of the vector field to a nonzero level surface is birationally equivalent to a linear vector field in~$\mathbf{CP}^2$}. The vector field has the quasihomogeneous polynomial first integral~$Q=x^3-3yx+z$. Let~$\Sigma=Q^{-1}(1)$. Its is biholomorphic to~$\mathbf{C}^2$  for it can be parametrized by~$(x,y)\mapsto (x,y,3xy-x^3)$. Let
$$\ell=(x^2-y+x+1,x^2-y+\rho^2 x+\rho,x^2-y+\rho  x+\rho^2),$$
for~$\rho$ a primitive cubic root of unity.
The function $j:\Sigma\to \mathbf{CP}^2$, given by the restriction of
$(x,y,z)\to[\ell_1:\ell_2:\ell_3]$, maps~$V|_\Sigma$ to the restriction to the image of the vector field~$\Delta$ given in coordinates~$[\xi:\zeta:1]$ by
$(\rho-1)\xi\partial/\partial \xi+(2\rho+1)\zeta\partial/\partial\zeta$. There is an action of~$\mathbf{Z}/3\mathbf{Z}$ on~$\mathbf{C}^3$ given by~$(x,y,z)\mapsto (\rho x,\rho^2 y, z)$ that preserves the level sets of~$Q$. There is also an action of~$\mathbf{Z}/3\mathbf{Z}$ on~$\mathbf{CP}^2$ given by~$[\xi:\zeta:1]\mapsto[\zeta:1:\xi]$. The map~$j$ is equivariant since~$\ell_k(\rho x,\rho^2 y)=\rho^2\ell_{k+1}(x,y)$.  The inverse of~$j$ (as a meromorphic map) is given by
$(x,y)=\left(f,\Delta\cdot f \right)$ for
$$f=\frac{(\rho+1)(\rho^2\zeta-\rho\xi+1)}{(\rho+\rho^2\zeta+\xi)}.$$ 
This proves the claim.

There is an invariant surface for~$V$ given by the zero locus of~$H=3y^2x^2-y^3-3xyz+z^2$ (this is, $H|V\cdot H$) which was used to construct the above solution. The intersection of~$\Sigma$ with~$H$ gives an invariant curve for the restriction of~$V$ to~$\Sigma$. This curve is not irreducible: $Q-1$ divides~$\ell_1\ell_2\ell_3-H$ within~$\mathbf{C}[x,y,z]$, this is, within~$\Sigma$, the zero locus of~$H$ agrees with the zero locus of~$\ell_1\ell_2\ell_3$. Within~$\Sigma$ (parametrized by~$x$ and~$y$),  $V|_\Sigma \cdot \ell_k|_\Sigma=(x-\rho^{1- k})\ell_k|_\Sigma$, equivalent to the above formulation (in particular, each one of the irreducible components of the intersection of~$\Sigma$ with~$H=0$ is mapped to a coordinate line  in~$\mathbf{CP}^2$).

\subsection*{Chazy V}
\emph{There is a quasihomogeneous polynomial first integral. A nonzero level surface embeds into~$\mathbf{CP}^1\times\mathbf{CP}^1$, while the restriction of the vector field is mapped to a product holomorphic vector field}. We have the quasihomogeneous first integral $Q=x^4-4x^2y+2zx-y^2$. Let~$\Sigma=Q^{-1}(1)$. Let~$\Delta$ be the
vector field on~$\mathbf{CP}^1\times\mathbf{CP}^1$ given in coordinates~$([\xi:1],[\zeta:1])$  by
$(i-1)(i\xi\partial/\partial\xi+\zeta\partial/\partial \zeta)$ and let
$$f=\frac{(\xi-1)(\zeta-1)}{\xi\zeta-i\xi-i\zeta+1}.$$
The mapping~$\Phi$ defined by~$([\xi:1],[\zeta:1])\dashrightarrow  (f,\Delta\cdot f,\Delta^2 \cdot f)$
maps $(\mathbf{CP}^1\times \mathbf{CP}^1,Z)$ to~$(\Sigma,V|_\Sigma)$. There is an action of~$\mathbf{Z}/4\mathbf{Z}$ on~$\mathbf{C}^3$ generated by~$(x,y,z)\mapsto (i x,- y, -i z)$ that preserves the level sets of~$Q$. There is also an action of~$\mathbf{Z}/4\mathbf{Z}$ on~$\mathbf{CP}^1\times\mathbf{CP}^1$ given by~$([\xi:1],[\zeta:1])\mapsto([\zeta:1],[1:\xi])$. The map~$\Phi$ is equivariant. 

The field~$V$ has an invariant irreducible hypersurface given by the vanishing of
$H=2y^2x^2-2xyz+z^2-2y^3$, that was used to construct the above solution: when parametrizing~$\Sigma$ by~$x$ and~$y$, $H|_\Sigma=\Pi_{k=1}^{4} C_i$ for
$$C_1=x^2+x+i+ix-y,\;C_2=x^2-ix-i+x-y,$$
$$C_3=x^2-x+i-ix-y,\;C_4=x^2+ix-i-x-y.$$
The inverse of~$\Phi$ (as a meromorphic map) is given by~$[C_1/C_3:C_2/C_4]$. The intersection of~$\{H=0\}$ with~$\Sigma$ has four irreducible components that correspond to the orbit (under the  action of~$\mathbf{Z}/4\mathbf{Z}$) of the line~$\{\xi=0\}\subset\mathbf{CP}\times\mathbf{CP}^1$.

\subsection*{Chazy VI}\emph{There is a quasihomogeneous polynomial first integral. The restriction of the vector field to a nonzero level surface is birationally equivalent to a linear vector field in~$\mathbf{CP}^2$}. The vector field~$V$ has the first integral $$Q=x^6-6x^4y+6zx^3-15x^2y^2+6xyz+8y^3-3z^2.$$
Let~$\Sigma=Q^{-1}(1)$ and
consider the birational mapping~$\Phi:\mathbf{CP}^2\dashrightarrow   \Sigma$ given by
$[\xi:\zeta:1] \to  (f,\Delta\cdot f,\Delta^2\cdot f)$ for
$\Delta=-\rho^2\xi\partial/\partial \xi+\zeta\partial/\partial \zeta,$
$$f=\frac{\xi(1-\zeta^2)+\rho  (\zeta^2-\xi^2)+\rho^2\zeta(\xi^2-1)}{\xi^2\zeta+\xi\zeta^2+\zeta^2 +\zeta +\xi +\xi^2 -6\xi\zeta},$$
for~$\rho$ a primitive cubic root of unity. We have~$\Phi(\mathbf{CP}^2,\Delta)=(\Sigma,V|_{\Sigma})$. The group~$\mathbf{Z}/6\mathbf{Z}$ acts upon~$\Sigma$ by the restriction of
the order six transformation~$(x,y,z)\mapsto(-\rho^2 x,\rho y,-z)$
and upon~$\mathbf{CP}^2$ by the order six birational transformation
$[\xi:\zeta:\nu]\mapsto[\xi\zeta:\nu\zeta:\xi\nu]$, that factors as a Cremona involution and a cyclic permutation of the variables. The mapping~$\Phi$ is equivariant with respect to these actions. It has poles along the denominator of $f$ and the points $[1:1:1]$, $[1:0:0]$, $[0:1:0]$ and~$[0:0:1]$ are indeterminacy ones. Let $S$ be the surface corresponding to the blow up of the last three.  It has a configuration of six rational curves, invariant by the vector field,  of self-intersection~$-1$ (a hexagon formed by the three coordinate axes plus the exceptional divisors of the blowup of their intersections). The group of order six acts holomorphically on~$S$ and transitively upon the curves in the hexagon. 

To prove that~$\Phi:\mathbf{CP}^2\to\Sigma$ is invertible as a meromorphic map, notice that~$\Phi$ is of the form $(P_1\Lambda^{-1}, P_2\Lambda^{-2}, P_3\Lambda^{-3})$ with~$P_i\in\mathbf{C}[\zeta,\xi]$ of degree $3i$ and $\Lambda$  a polynomial of degree 3. The invertibility of~$\Phi$ is equivalent to the fact that there are six points in~$\mathbf{CP}^2$ whose image under~$\Phi$ intersect the curve~$(t,\alpha t^2,\beta t^3)$ for generic (but fixed) $\alpha$ and~$\beta$. These points of intersection belong to the common zeroes of the polynomial equations~$E_2=P_2-\alpha P_1^2$ and~$E_3=P_3-\beta P_1^3$. Counted with multiplicity, there are~$54$ solutions to these equations in~$\mathbf{CP}^2$, among which we have:
\begin{itemize}
 \item The point~$[0:0:1]$ with multiplicity~$6$ ($E_2$ has two smooth branches at this point and $E_3$ has three; these five branches have different tangents). By symmetry, the points~$[0:1:0]$ and~$[1:0:0]$ have also multiplicity~$6$ (these account for~$18$ solutions).
 \item The point~$[1:1:1]$ with multiplicity~$30$. The polynomial $E_2$ has two smooth (tangent) branches at this point. After calculating the first terms in the Puiseux parametrization of each branch and substituting in~$E_3$, we obtain that the multiplicity of the intersection of each branch of~$E_2=0$ with~$E_3=0$ is at least~$15$.
\end{itemize}
This gives us at least~$18+30=48$ common solutions to~$E_2$ and~$E_3$ (all of them indeterminacy points of~$\Phi$) and leaves at most~$54-48=6$ other points. We conclude that~$\Phi$ is invertible.

The vector field~$V$ has an invariant hypersurface given by the zero locus of~$H=y^2x^2-3y^3-xyz+z^2$, which was used to find the above solution. The intersection of this surface with~$\Sigma$ is the image of the hexagon under~$\Phi$.

\subsection*{Chazy VII}  \emph{The vector field is completely integrable. The generic orbit is an elliptic curve.} We have the quasihomogenous  first integrals~$g_2=\frac{4}{3}(x^4+2x^2y+y^2-2xz)$ and
$g_3=-\frac{4}{27}(2x^6+6x^4y+6y^2x^2-2y^3-6x^3z-6zxy+3z^2)$. If we set~$P=\frac{1}{3}x^2+\frac{2}{3}y$, we have the Weierstrass differential equation~$(P')^2=4P^3-g_2P-g_3$ and thus the generic fiber is an elliptic curve.

\subsection*{Chazy VIII} \emph{The vector field is completely integrable. The generic orbit is an elliptic curve.} We have the quasihomogenous first integrals~$a=z-2x^3$ and~$b=y^2-2zx+3x^4$. We have
$y^2=x^4+2(z-2x^3)x+(y^2-2zx+3x^4)$, equivalent to~$(\phi')^2=\phi^4+2a\phi+b$. The curves are elliptic.

\subsection*{Chazy IX}\emph{There is a quasihomogenous first integral. The generic level surface embeds into the Jacobian of the curve~$\xi^2=\zeta^5-1$ and, under this embedding, the vector field becomes one of the holomorphic vector fields in the Jacobian that is preserved, up to a constant factor, by the action of the  automorphism group of the curve}. Part of what follows is a rephrasing of Cosgrove~\cite[\S 3]{cosgrove-2}. Consider the genus two curve~$C$ given by
$\xi^2=\zeta^5-1$.
It has a cyclic group of order~$10$ of automorphisms generated by
\begin{equation}\label{sym10}(\zeta,\xi)\mapsto (\omega \zeta,-\xi),\end{equation}
for~$\omega$ a primitive fifth root of unity.
We will follow Mumford's account of Jacobi's work for the algebraic construction of~$\mathcal{J}$, the Jacobian of~$C$~\cite{theta2}. Consider~$\mathbf{C}^4$ with coordinates~$u_1$, $u_2$, $v_1$ and~$v_2$. The complement of the theta divisor within $\mathcal{J}$ may be embedded in~$\mathbf{C}^4$ and its image is the variety~$Z$ given by the equations
$$-1-2u_2^2u_1+u_2u_1^3+u_2v_1^2-v_2^2=0,\;-2v_1v_2+u_2^2-3u_2u_1^2+u_1^4+u_1v_1^2=0.$$
These equations give  the necessary and sufficient conditions for $t^2+u_1t+u_2$ to divide~$(t^5-1)-(v_1t+v_2)^2$. The order~$10$ transformation induced by~(\ref{sym10}) is given by
$$\label{a5symmetries}(u_1,u_2,v_1,v_2)\mapsto(\omega u_1,\omega^2 u_2,-\omega^4 v_1,-v_2).$$
The Jacobian of~$C$ is the only Abelian surface having a cyclic group of order~$10$ of automorphisms~\cite{fujiki}.
The vector fields in~$\mathcal{J}$ are generated by (extending) the commuting vector fields on~$Z$
$$\Delta_1=v_1\del{u_1}+v_2\del{u_2}+(u_2-{\textstyle \frac{3}{2}u_1^2)}\del{v_1}+({\textstyle\frac{1}{2}}v_1^2+{\textstyle \frac{1}{2}}u_1^3-2u_1u_2)\del{v_2},$$
$$\Delta_2=v_2\del{u_1}+(v_2u_1-u_2v_1)\del{u_2}+({\textstyle \frac{1}{2}} v_1^2+{\textstyle \frac{1}{2}}u_1^3-2u_1u_2)\del{v_1}+{\textstyle \frac{1}{2}}(u_1^4-u_2u_1^2+u_1v_1^2-2u_2^2)\del{v_2}.$$
Both vector fields are preserved, up to a constant factor, by the transformation~(\ref{a5symmetries}).\\

Let us now come back to Chazy~XI. The vector field~$V$ is divergence-free and has the quasihomogeneous first integral (of degree~10)
$$Q=-10 x z^3+5(36 x^4+y^2) z^2+60(27 x^4+12 y x^2+ 2y^2) xyz+2916 x^{10}-48 (3 x^2+y)^5.$$
Let~$f$ be the meromorphic function in~$Z$ given by
$$f=\frac{1}{6}\frac{(\sqrt{5}-1)(u_1v_1-2v_2)}{2u_1^2+(\sqrt{5}-3)u_2}.$$
Let~$\Phi:Z\dashrightarrow  \mathbf{C}^3$ be given by~$\Phi(u,v)=(f,\Delta_1 \cdot f,\Delta_1^2 \cdot f)$.
Its image belongs to a non-zero level surface~$\Sigma$ of~$Q$ (depending upon the determination of~$\sqrt{5}$) and maps~$\Delta_1$ to the restriction to~$\Sigma$ of~$V$. It is equivariant with respect to the group of transformations generated by
\begin{equation}\label{massym}(x,y,z)\mapsto (-\omega x, \omega^2 y,-\omega^3 z),\end{equation}
acting on~$\Sigma$, and the group~(\ref{a5symmetries}) acting on~$Z$. The inverse of~$\Phi$, $\Phi^{-1}:\Sigma\dashrightarrow Z$ is given by the restriction to~$\Sigma$ of the mapping
\begin{multline*}(x,y,z)\mapsto {\textstyle \frac{1}{5}} (-3[5+\sqrt{5}][6x^2+(1+\sqrt{5})y],
18[(27+4\sqrt{5})y^2+162x^4+54x^2y-2(\sqrt{5}+5)yz],
\\-3[5+\sqrt{5}][(1+\sqrt{5})z+12xy],
\\-{\textstyle \frac{36}{5}}[5+\sqrt{5}][270x^5+36(5+\sqrt{5})x^3y+(3\sqrt{5}-15)zx^2+6(5+\sqrt{5})y^2x-(5
+2\sqrt{5})zy]),\end{multline*}
and~$\Phi$ is thus a bimeromorphism. The image by~$\Phi$ of~$(\frac{5}{12}-\frac{1}{4}\sqrt{5})\Delta_2$ is
\begin{multline*}W=(36 x^4+12 x^2 y+2 y^2-3 x z)\del{x}+(yz-162x^5-72x^3y-12y^2x+12zx^2)\del{y}+\\+(648 x^6+108 x^4 y-72 x^3 z+z^2)\del{z},\end{multline*}
which commutes with~$V$ and satisfies~$[L,W]=3W$ (is quasihomogeneous of degree~4). It is divergence-free, has a first integral in~$Q$ and is also preserved, up to a constant factor, by~(\ref{massym}).

In order to give explicit functions of~$t$ solving the equation, we would have to parametrize the orbits of~$\Delta_1$. A priori, this can be done through theta functions on~$\mathcal{J}$~\cite{theta2} (there is a  polarization in~$\mathcal{J}$ invariant by the automorphism group of~$C$~\cite{fujiki}). We will not pursue this approach.

%%%%%%%%%%%%%%%%%%%%%%%%%%%%%%%%%%

\subsection*{Chazy X} \emph{There is a quasihomogenous first integral. The generic level surface identifies to an open subset in the square of the elliptic curve admitting an order four multiplication. This identification maps the Chazy vector field to a holomorphic one that is preserved,  up to a constant factor, by a cyclic group of automorphisms of the surface.} Part of our development rephrases the work of Cosgrove~\cite[\S 4]{cosgrove-2}. Let~$u(z)=\mathrm{sn}(z,i)$, the Jacobi elliptic function that satisfies the differential equation
$(u')^2=1-u^4$. The period lattice~$\Lambda$ of~$u$ is generated by~$4K$ and~$2iK'$ for~$K=\int_0^1(\sqrt{1-\zeta^4})^{-1}d\zeta\in\mathbf{R}$ and~$K'=K-iK$. For the half-periods we have
\begin{equation}\label{halfperiod}u(z+2K)=-u(z),\; u(z+iK')=-\frac{i}{u(z)}, \;  u(z+2K+iK')=\frac{i}{u(z)}.\end{equation}
Moreover,
\begin{equation}\label{orderfour}u(i\xi)=iu(\xi),\; u'(i\xi)=u'(\xi).\end{equation}

The vector field~$V$ is divergence-free and admits the quasihomogeneous first integral of the twelfth degree
\begin{multline*} Q=(94392+52164  \sqrt{3}) x^{12}+(263088  \sqrt{3}+423792) x^{10} y+(322704  \sqrt{3}+617544) x^8 y^2+\\+(388584+254712  \sqrt{3}) x^6 y^3+(66492  \sqrt{3}+143136) x^4 y^4+(5688  \sqrt{3}+3240) x^2 y^5+\\+(-8480-4992  \sqrt{3}) y^6+(-127512  \sqrt{3}-163944) x^7 y z+(-77616  \sqrt{3}-221760) x^5 y^2 z+\\+(-56232-63096  \sqrt{3}) x^3 y^3 z+(-10032+3168  \sqrt{3}) x y^4 z+(-16632  \sqrt{3}-47520) x^6 z^2+\\+(15444  \sqrt{3}+1188) x^4 y z^2+(20064-6336  \sqrt{3}) x^2 y^2 z^2+(396+5148  \sqrt{3}) y^3 z^2+22264 x^3 z^3+\\+(1815  \sqrt{3}-4356) z^4.
\end{multline*}

Let~$\rho_1=\frac{1}{2}(1+i+ \sqrt{3}+i \sqrt{3})$ and $\rho_2=\frac{1}{2}(1+i- \sqrt{3}-i \sqrt{3})$ and consider, in~$\mathbf{C}^2$, the commuting vector fields $$\Delta_1=\rho_1\del{\xi}+\del{\zeta}, \; \Delta_2=\rho_2\del{\xi}+\del{\zeta}$$
and the meromorphic function
\begin{multline*}f(\xi,\zeta)=\frac{ u'(\zeta)-i(2+ \sqrt{3})u'(\xi)}{u(\zeta)-\rho_1 u(\xi)}+\\+\frac{\rho_1[i u(\xi) u'(\xi)+(i+1) u(\zeta) u'(\xi)+u(\xi) u'(\zeta)][i u(\xi)+\rho_2 u(\zeta)] u(\xi)}{2-u(\xi)^2[u(\xi)^2+(1-i) u(\xi)u(\zeta)+i  u(\zeta)^2]}.\end{multline*}
The function~$\Phi:\mathbf{C}^2\dashrightarrow  \mathbf{C}^3$ given by~$\Phi(\xi,\zeta)=(f, \Delta_1\cdot f,\Delta_1^2\cdot f)$ takes values on a non-zero level surface~$\Sigma$ of~$Q$ and maps~$\Delta_1$ to~$V$. Under~$\Phi$, the vector field~$(13904+8008\sqrt{3})\Delta_2$ is mapped to the vector field
\begin{multline*}W=[(540+288  \sqrt{3}) x^6+(432+270  \sqrt{3}) x^4 y+(168+94  \sqrt{3}) x^2 y^2-(32+20  \sqrt{3}) y^3+ \\ +33 z^2-22(5 \sqrt{3}+12) x^3 z-44(\sqrt{3}+1) x y z] \del{x} +
[(-774  \sqrt{3}-1278) x^7+\\(594+396  \sqrt{3}) x^4 z -(1260+672  \sqrt{3}) x^5 y-(534+342  \sqrt{3}) x^3 y^2-(4  \sqrt{3}+24) x y^3+\\ +22(5 \sqrt{3}+12) x^2 y z+22(  \sqrt{3}+1) y^2 z-44(  \sqrt{3}+1) x z^2]\del{y} +[162(13  \sqrt{3}+23) x^8+\\ +108(31+18  \sqrt{3}) x^6 y+6(175 \sqrt{3}+291) x^4 y^2+24(7+3  \sqrt{3}) x^2 y^3+(156+92  \sqrt{3}) y^4+\\ +22(5 \sqrt{3}+12) x^2 z^2-132(15+8  \sqrt{3}) x^5 z-132(5+3  \sqrt{3}) x^3 y z-(176  \sqrt{3}+264) x y^2 z]\del{z},\end{multline*}
which is divergence-free, commutes with~$V$, has in~$Q$ a first integral and satisfies~$[L,W]=5W$ (it is quasihomogeneous of degree~6).

The function $f$ (together with its derivatives along~$\Delta_1$) is invariant with respect to the lattice~$\Lambda^2\subset\mathbf{C}^2$ generated, in each factor, by the periods of~$u$. By the relation~(\ref{halfperiod}), $f$~has also the period~$\tau=(2K,2K)\notin \Lambda^2$ (this period is absent in Cosgrove's account). Let~$\Gamma\subset\mathbf{C}^2$ be the lattice generated by~$\Lambda^2$ and~$\tau$ (it is still a lattice since~$2\tau\in\Lambda^2$). The mapping~$\Phi$ induces a mapping~$\Phi^\Gamma:\mathbf{C}^2/\Gamma\dashrightarrow\Sigma$.

By the quasihomogeneity of the first integral, any level surface of~$Q$ is mapped to itself under the order~$12$ transformation generated by
\begin{equation}\label{symor12}(x,y,z)\mapsto (\rho x, \rho^2 y,\rho^3 z)\end{equation}
for~$\rho$ a primitive~$12^\text{th}$ root of unity. This transformation multiplies the vector field by a suitable power of~$\rho$. We claim that~$\Phi$ is equivariant with respect to the action of~$\mathbf{Z}/12\mathbf{Z}$ on~$\Sigma$ given by the above formula and the action on~$\mathbf{C}^2$ generated by the linear transformations
$$A=\left(\begin{array}{cc} i & 0 \\ 0  & i\end{array}\right), \; B=\left(\begin{array}{rr}\frac{1}{2}[i-1] & \frac{1}{2}[i-1] \\  \frac{1}{2}[i+1] & -\frac{1}{2}[i+1]\end{array}\right),$$
that generate a cyclic group of order~$12$ ($A^4=1$, $B^3=1$, $AB=BA$). The assertion concerning~$A$ follows in a straightforward way: by~(\ref{orderfour}), $f(i\xi,i\zeta)=-if(\xi,\zeta)$. The assertion concerning~$B$ follows by the fact that~$B$ and~$B^2$ are the unique order three transformations preserving~$\Delta_1$ and~$\Delta_2$ up to a constant (and different) factor and hence correspond, in~$\Sigma$, to the order three transformations induced by~(\ref{symor12}). This also follows, in a more complicated way, by the addition formul{\ae} for Jacobi elliptic functions.

Notice that the action of~$B$  preserves~$\Gamma$ (but does not preserve~$\Lambda^2$). Hence, $\Phi^{\Gamma}:\mathbf{C}^2/\Gamma\dashrightarrow \Sigma$ is equivariant with respect to the action of~$\mathbf{Z}/12\mathbf{Z}$.
The lattice~$\Gamma$ may be generated by~$e_1$, $Ae_1$, $e_2$ and~$Ae_2$ for $e_1=2K(0,1+i)$ and $e_2=2K(1,i)$. Hence~$\mathbf{C}^2/\Gamma$ is isomorphic to the square of the elliptic curve admitting an order four multiplication. In the basis of~$\mathbf{C}^2$ given by~$e_1$ and~$e_2$ the order three transformation given by~$B$ is given by
$\left(\begin{array}{rr} 0 & 1 \\ -1 & -1 \end{array}\right)$, agreeing with the conventions given by Fujiki~\cite[Table~6]{fujiki} for this Abelian surface with cyclic automorphism group.\\

It remains to be shown that~$\Phi^{\Gamma}:\mathbf{C}^2/\Gamma\dashrightarrow\Sigma$ is a bimeromorphism, this is, to prove that if~$\Gamma'\subset\mathbf{C}^2$ is a lattice such that~$\Gamma\subset\Gamma'$ and such that~$\Phi^{\Gamma'}:\mathbf{C}^2/\Gamma'\dashrightarrow\Sigma$ is well-defined, $\Gamma'=\Gamma$. Let~$\Upsilon:\mathbf{C}^2/\Lambda^2\to\mathbf{C}^2/\Gamma'$ be the natural quotient and let, for~$i=1,2$, $\mathcal{C}_i=\mathbf{C}/\Lambda$ so that~$\mathbf{C}^2/\Lambda^2=\mathcal{C}_1\times \mathcal{C}_2$. If~$p\in\Gamma'$ then~$p=p_1+p_2$ for~$p_i\in\mathcal{C}_i$. In~$\mathbf{C}^2/\Gamma'$, $\Upsilon(p_1)+\Upsilon(p_2)=\Upsilon(p)=0$ and thus~$\Upsilon(p_1)=\Upsilon(-p_2)$. Reciprocally, if~$\Upsilon(p_1)+\Upsilon(p_2)=0$ then~$p_1-p_2$ belongs to~$\Gamma'$. Hence, in order to understand the elements of~$\Gamma'$ it is sufficient to consider the couples~$(p_1,p_2)\in\mathcal{C}_1\times \mathcal{C}_2$ such that~$\Upsilon(p_1)=\Upsilon(-p_2)$.

Let~$\mathbf{CP}^2_L$ be the weighted projective space associated to~$L$ and consider the natural mapping~$\Pi:\mathbf{C}^3\setminus\{0\}\to \mathbf{CP}^2_L$. 
We will write~$[x:y:z]_L$ for~$\Pi(x,y,z)$.  We have a mapping~$\Pi\circ \Phi^{\Lambda^2}:\mathbf{C}^2/\Lambda^2\to\mathbf{CP}^2_L$. We can factor $\Phi^{\Lambda^2}$ as~$\Phi^{\Gamma'}\circ\Upsilon$. In particular, if~$\Upsilon(p_1)=\Upsilon(-p_2)$, $\Pi\circ\Phi^{\Lambda^2}(p_1)=\Pi\circ\Phi^{\Lambda^2}(-p_2)$.

The curve~$\Pi\circ\Phi^{\Lambda^2}(0,\zeta)$, image of~$\mathcal{C}_2$, is a curve of degree~$9$ parametrized by
\begin{multline}\label{parc2}t\mapsto\left[3t-i\sqrt{3}:{\textstyle \frac{3}{22}}(i\sqrt{3}+5)(11t^3-i(14 \sqrt{3} -4)t^2 -(17-10\sqrt{3})t+12i\sqrt{3}+6i):\right. \\ \left.3i(2\sqrt{3}+3)(-11t^4+(17\sqrt{3}-26)it^3+(57-16\sqrt{3})t^2+(26-17\sqrt{3})it-14)\right]_L \end{multline}
for~$t=u'(\zeta)$. On the other hand, the curve~$H\circ\Phi^{\Lambda^2}(\xi,0)$, image of~$\mathcal{C}_1$, has the quasihomogeneous equation
\begin{equation}\label{eqc1}(3+6\sqrt{3})x^6-11xyz+(42+18\sqrt{3})x^4y+11x^3z+(6\sqrt{3}+14)y^3-{\textstyle \frac{11}{2}}x^2y^2+(31+18\sqrt{3})x^2y^2\end{equation} and is parametrized by
\begin{multline}\label{parc1}s\mapsto\left[s-2i-i\sqrt{3}:{\textstyle\frac{1}{6}}(\sqrt{3}+3)(3s^2+2i[3+\sqrt{3}]s+6\sqrt{3}+3):\right.\\ \left.(\sqrt{3}+2)(3s^3-i\sqrt{3}s^2-[8\sqrt{3}+11]s+9i\sqrt(3)+8i)\right]_L, \end{multline}
with~$s=u'(\xi)$. By evaluating the equation~(\ref{eqc1}) on the parametrization~(\ref{parc2}), we obtain, up to a constant multiple,
$(t-1)^2(t+1)^2(t+i)^5$. The curves intersect only at the roots of this polynomial (at~$t=\infty$ we attain the point~$[0:1:0]_L$, which does not satisfy the equation).
\begin{itemize}
\item For~$t\in\{-1,1\}$ we attain the point $p=[1:-1:2]_L$, corresponding to~$s\in\{-1,1\}$. The point~$p$ is the image of the four points in~$\mathbf{C}^2/\Lambda^2$ for which~$u(\xi)=0$ and~$u(\zeta)=0$. The function~$u$ has degree~$2$ and it vanishes at~$0$ and~$2K$ (inequivalent mod~$\Lambda$), so we have the possible periods~$(0,0)$, $(2K,2K)$ (which are indeed periods belonging to~$\Gamma$) and $(0,2K)$ and~$(2K,0)$. If~$(0,2K)$ were a period then, by~(\ref{halfperiod}), the parametrization~(\ref{parc2}) would be invariant under a change of sign of~$t$ (analogously, if~$(2K,0)$ were a period, (\ref{parc1}) would be invariant by a change of sign of~$s$). We must conclude that no extra periods arise for~$t=\pm 1$.
\item For~$t=-i$, we attain the point~$[1:\frac{3}{2}+\frac{1}{2}\sqrt{3}:6+3\sqrt{3}]_L$, corresponding to~$s=\infty$. The two poles of~$u'$, correspond, by~(\ref{halfperiod}), to~$iK'$ and~$2K+iK'$ (half-periods of~$\Lambda$). The possible periods are of the form~$(\nu,\mu)$, where~$u'(\nu)=-i$ (and thus~$u^4(\nu)=2$) and~$\mu$ is one of the half-periods of~$\Lambda$ where~$u$ has poles. If some of these were periods, their doubles, of the form~$(2\nu,0)$, should be periods as well, this is, $2\nu\in\Lambda$. But then, by~(\ref{halfperiod}), $u(\nu)\in\{0,\infty\}$.
\end{itemize}
We conclude that~$\Gamma'=\Gamma$ and thus that~$\Phi:\mathbf{C}^2/\Gamma\dashrightarrow\Sigma$ is a bimeromorphism.

\subsection*{Chazy XI} \emph{There is a first integral. The generic level surface compactifies as a rational fibration over an elliptic curve, where the vector field extends as a holomorphic one.}  We will suppose that~$k>0$. Let~$\wp$ be the Weierstrass function satisfying~$(\wp')^2=4\wp^3-1$, $\Lambda\subset\mathbf{C}$ its lattice of periods and let~$E=\mathbf{C}/\Lambda$. Consider~$\mathbf{CP}^1\times E$ as the trivial fibration~$\Pi:\mathbf{CP}^1\times E\to E$. Give to~$\mathbf{CP}^1\times E$ coordinates~$([\phi:1], t \mod \Lambda)$. Consider,  on~$\mathbf{CP}^1\times E$, the meromorphic  vector field
\begin{equation}\label{delta11}\Delta={\textstyle\frac{1}{2}}([1-k]\phi^2+[1+k]\wp(t))\del{\phi}+\del{t}.\end{equation}
It is transverse to the fibers of~$\Pi$, except for the fiber~$\Pi^{-1}(0)$, where it has poles. Under~$\Pi$, $\Delta$ is mapped to the vector field $\partial/\partial t$ on~$E$. \\

Under the quasihomogeneous change of variables $$(X,Y,Z)=({\textstyle\frac{1}{2}}[k+1]x,y+{\textstyle\frac{1}{4}}[k^2-1]x^2,z+{\textstyle\frac{1}{2}}[k^2-1]xy), $$
the vector field~$V$ becomes
$${\textstyle\frac{1}{2}}([1-k]X^2+[1+k]Y)\del{X}+Z\del{Y}+6Y^2\del{Z}.$$ The latter has the first integral~$g_3=4Y^3-Z^2$.   Let~$\Sigma=g_3^{-1}(1)$.
Let~$\Phi:\mathbf{CP}^1\times E\to\Sigma$ be given by~$\Phi([\phi:1], t)=(\phi,\wp(t),\wp'(t))$. This map is well-defined and maps~$\Delta$ to the restriction of~$V$ to~$\Sigma$. In other words, $(\Sigma,V|_\Sigma)$ embeds into~$(\mathbf{CP}^1\times E,\Delta)$.\\

In this way, $V$ is semicomplete if and only if the vector field~$\Delta$ on~$\mathbf{CP}^1\times E$ is semicomplete. We will prove this by showing that, when~$k\notin 6\mathbf{Z}$, $\Pi:\mathbf{CP}^1\times E\to E$ can be birationally transformed into a (no longer trivial) rational fibration  where~$\Delta$ becomes holomorphic (any holomorphic vector field on a compact manifold is automatically complete).\\

In order to understand~$\Delta$ in the neighborhood of the fiber~$t=0$ (where~$\Delta$ has poles), notice that, after the change of variables~$\xi=t\phi+1$ (that amounts to birationally modifying the fibration), the vector field~(\ref{delta11}) becomes
\begin{equation}\label{ricchazy11} \left \{k\xi+{\textstyle\frac{1}{2}}[1-k]\xi^2+{\textstyle\frac{1}{2}}[1+k](t^2\wp(t)-1)\right\}t^{-1}\del{\xi}+\del{t}.\end{equation}
Notice that, in the above equation, $f(t)=t^2\wp(t)-1=\frac{1}{28}t^6+\cdots$ is a holomorphic function vanishing at the origin. The above equation has a pole at~$t=0$ that we will try to eliminate by transforming the fibration (to illustrate what follows, if we set~$\xi=\mu t^k$, the vector field  becomes
$${\textstyle\frac{1}{2}}\left \{[1-k]\mu^2t^{k-1}+[1+k]t^{-k-1}f(t)\right\}\del{\mu}+\del{t},$$
which is holomorphic and transverse to~$t=0$ if~$k<6$). The foliation induced by~(\ref{ricchazy11}) has a singularity at the origin that is in the Poincar\'e domain and is resonant. Let~$F$ denote the fiber~$t=0$. It has self-intersection~$0$. Blowing-up the origin produces an exceptional divisor~$D$ of self-intersection~$-1$ and the strict transform of~$F$ becomes a curve of self-intersection~$-1$ which can be, on its turn, blown down, making~$D$ a curve of self-intersection~$0$ (this transformation is known as a \emph{flip} on the fiber).
By the Poincar\'e-Dulac Theorem, there exist coordinates~$(x,y)$ in a neighborhood of the origin where the vector field~(\ref{delta11}) has the expression
\begin{equation}\label{poincaredulac}f(x,y)\left[\del{x}+x^{-1}(ky+\delta x^k)\del{y}\right]\end{equation}
for some holomorphic function~$f$ that does not vanish at the origin and some~$\delta\in\{0,1\}$.  In the chart given by~$y=\sigma x$ of the blowup  of the origin, the last vector field becomes~$f(x,\sigma x)[\partial/\partial{x}+x^{-1}([k-1]\sigma+\delta x^{k-1})\partial/\partial{\sigma}]$. Thus, by successively flipping the fiber, we get to the case~$k=1$. In this case, after blowing up, if~$\delta=0$, the meromorphic vector field  becomes the holomorphic one~$f(x,\sigma x)\partial/\partial{x}$, transverse to the exceptional divisor~$\{x=0\}$. Hence, in order to prove that the vector field~(\ref{delta11}) is semicomplete for~$k\notin6\mathbf{Z}$, it is sufficient to prove that in the corresponding Poincar\'e-Dulac normal form for the origin, the value of~$\delta$ is zero. This is equivalent to the existence of a formal solution of~(\ref{delta11}) \emph{up to order~$k$}. In fact, if we propose a formal solution~$(x(t),y(t))=(t,\sum_{i=1}^\infty a_it^i)$ we get (supposing without loss of generality that~$f\equiv 1$) the equations~$ia_i=ka_i$ for~$k\neq i$ (hence~$a_i=0$ if~$i\neq k$) and~$ka_k=ka_k+\delta$: there is a formal solution up to order~$k$ if and only if~$\delta=0$.

If  we set~$\mu(t^6)=\xi(t)$ in equation~(\ref{ricchazy11}), it becomes
$$6t\frac{d\mu}{d t}=k\mu+{\textstyle\frac{1}{2}}[1-k]\mu^2+{\textstyle\frac{1}{2}}[1+k]\rho(t),$$
for the holomorphic function~$\rho$ such that~$\rho(t^6)=t^2\wp(t)-1$.
If~$k\notin 6\mathbf{Z}$, there is no obstruction for the existence of a formal solution, for~$k/6$ is not a positive integer (see~\cite[\S~12.6]{ince} for details) and hence there is a formal solution for~$\mu(t)$, yielding a formal solution for~$\phi(t)$. Under~$k$ successive flips of the fiber the vector field becomes holomorphic and is transverse to the rational fibration: the vector field is semicomplete. By studying the formal solutions of this equation we should also be able to prove that there is no formal solution in the case where~$k\in 6\mathbf{Z}$, which would imply that not every solution of Chazy~XI for~$k\in 6\mathbf{Z}$ is single-valued. We are not aware of any published proof of this fact. The original work of Chazy~\cite{chazy}, which most authors refer to on this subject, does not make such a claim.

\subsection*{Chazy XII} In all these equations the corresponding vector field is semicomplete and from a solution~$\phi(t)$ and~$\left(\begin{array}{cc}a & b \\ c & d \end{array}\right)\in\mathrm{SL}(2,\mathbf{C})$, we obtain a new solution given by
$$\widetilde{\phi}(t)=\frac{1}{(ct+d)^2}\phi\left(\frac{at+b}{ct+d}\right)-\frac{6c}{ct+d}.$$
In this way we obtain the general solution. If~$k\in\{2,3,4,5\}$, every solution is a rational one and the system is completely integrable~\cite{chazy}: particular solutions are given by the logarithmic derivatives of~$P^{\frac{1}{2}k-3}$ for~$P$ equal to
\begin{eqnarray*}
t^2+1, & \text{if} & k=2 \\
3t^3+3t+1, & \text{if} & k=3 \\
t^5+t, & \text{if} & k=4 \\
t^{11}+11t^6-t, & \text{if} &  k=5
\end{eqnarray*}
and by the above method one may obtain the general solution. The interesting case corresponds thus to~$k>7$. In all these cases, the solutions are defined in a subset of~$\mathbf{C}$ that has a natural boundary (in particular, uncountable complement) so, according to~\cite[Corollary~D]{gui-reb}, there is no meromorphic first integral (and, in consequence, do not fall within the scope of Theorem~\ref{guillot-rebelo}).  Many works have been devoted to these equations. We will briefly state some of our own results in~\cite{guillot-sl2}.\\

Within the group of isometries of the hyperbolic plane generated by reflections on the sides of a triangle with internal angles $\pi/2$, $\pi/3$, $\pi/k$, let $T(2,3,k)\subset \mathrm{PSL}(2,\mathbf{R})$ be the (discrete) index two subgroup of orientation preserving ones. There is an embedding of $\mathbf{C}^3$ into a particular complex threefold~$M$ (which can be chosen compact if~$k\neq\infty$, and will be non-K\"ahler and of zero algebraic dimension) where~$V$ extends holomorphically to a \emph{complete} vector field (the complement of the image of~$\mathbf{C}^3$ in~$M$ has non-empty interior). This manifold~$M$ is endowed with an action of~$\mathrm{PSL}(2,\mathbf{C})$ and the stabilizer of its unique three-dimensional orbit is~$T(2,3,k)$. The trace in~$\mathbf{C}^3$ of the vector fields that generate the action in~$M$  are~$V$, $L$ and the vector field
$A=6\partial/\partial{x}+2x\partial/\partial{y}+6y\partial/\partial{z}$,
that satisfies the relation~$[A,V]=2L$: these three vector fields in~$\mathbf{C}^3$ generate a Lie algebra isomorphic to~$\mathfrak{sl}(2,\mathbf{C}$).

\section{Recovering Chazy's classification}\label{cuentas}

We will now give a proof of the fact that the only equations of the form~(\ref{eq:chazy}) that yield semicomplete vector fields are those in Chazy's list. This follows from exhibiting obstructions for semicompleteness in all the other equations. There are many more or less equivalent formulations (although different interpretations) of these obstructions, the first ones considered by Chazy himself. We will follow our own account in~\cite{guillot-semi} and  refer the reader to this article for details and proofs. We will consider the vector field~$V$ from equation~(\ref{field:chazy}). The key point is that  the vector field~(\ref{ell}) has $2i\pi$-periodic solutions and that the relation~$[L,V]=V$ holds. We have the following result~\cite[Corollary~2.6]{guillot-semi}:

\begin{prop}\label{crit} Let~$p\in\mathbf{C}^3$ be a point where~$L$ and~$V$ are collinear and where~$V$ does not vanish. Let~$c\in \mathbf{C}$ be such that~$L+cV$ vanishes at~$p$. If~$V$ is semicomplete then~$L+cV$ is linearizable around~$p$, its linear part is diagonalizable and its eigenvalues belong to~$\mathbf{Z}$. In particular, if one of these eigenvalues vanishes then
the locus of collinearity of~$L$ and~$V$ is at least two-dimensional. \end{prop}

The proof of this Proposition is based on the fact that if~$V$ is semicomplete then~$L+cV$ will not only be semicomplete but its solutions will be, like those of~$L$, $2i\pi$-periodic (in fact, $V$ and~$L$ generate a single-valued, although not complete, action of the affine group). We will now apply this criterion to the vector fields~(\ref{field:chazy}) in order to prove that, within this class, only Chazy's vector fields are semicomplete.\\

Let $\kappa_1$, $\kappa_2$ and~$\kappa_3$ be the roots of the polynomial
$$6\kappa^3-(2\delta+a_1)\kappa^2-a_2\kappa-a_3.$$
The vector fields~$V$ and~$L$ are linearly  dependent along the integral curve of~$L$ passing through the points $p_i=(1,\kappa_i,2\kappa_i^2)$. The field~$V$ will vanish along this orbit if and only if~$\kappa=0$. The vector fields~$V$ and~$L$ are linearly independent away from these three curves. Since this set is one-dimensional, the eigenvalues of Proposition~\ref{crit} relative to~$p_i$ (defined if~$\kappa_i\neq 0$), are non-zero in the semicomplete case. We have
\begin{eqnarray*}
 a_1 & = & 6(\kappa_1+\kappa_2+\kappa_3)-2\delta\\
a_2 & = &  -6(\kappa_1\kappa_2+\kappa_2\kappa_3+\kappa_3\kappa_1)\\
a_3 & = & 6 \kappa_1\kappa_2\kappa_3.
\end{eqnarray*}
If~$\phi(t)$ is a solution to~(\ref{eq:chazy}) then~$\phi(\mu^{-1} t)$ is a solution to the equation
$$\phi'''=\mu^3a_3\phi^4+\mu^2 a_2\phi^2\phi'+\mu a_1(\phi')^2+\mu\delta \phi\phi''$$
and in this way, the family of equations~(\ref{eq:chazy}) is actually a three-parameter one:
the vector field depends essentially on the ratio~$[\kappa_1:\kappa_2:\kappa_3:\delta]$. If~$\kappa_1\neq 0$ then, at the point~$p_1$, $L-\kappa_1^{-1}V$ will vanish. Its linear part has three eigenvalues, one of which is~$-1$. Let~$u_1$ and~$v_1$ be the other ones. Then~$t_1=u_1+v_1$ and~$d_1=u_1v_1$ are given by
\begin{equation}\label{defeigenvals}t_1=7-\frac{\delta}{\kappa_i},\; d_1=\frac{6(\kappa_1-\kappa_2)(\kappa_1-\kappa_3)}{\kappa_1^2}.\end{equation}

The six numbers~$t_i,d_i$ are functions of three parameters and in consequence, we must have three relations among them. A straightforward calculation shows that
\begin{equation}\label{relations}\sum_{i=1}^3\frac{1}{d_i}=\frac{1}{6},\; \; \sum_{i=1}^3\frac{t_i}{d_i}=\frac{7}{6}, \;\; \sum_{i=1}^3\frac{t_i^2}{d_i}=\frac{49}{6}.\end{equation}

These relations are analogue to the ones that exist for homogeneous vector fields in~$\mathbf{C}^n$ (see
\cite{guillot-fixe}). Notice that if~$t_1=7$ then~$t_i=7$ for every other~$i$ where~$\kappa_i\neq 0$. In the other cases, $t_i$ determines~$\kappa_i/\delta$ and thus, from a solution to~(\ref{relations}), we obtain, at most, a unique differential equation of the form~(\ref{eq:chazy}). Four cases arise:

\subsection{Case~1, $\kappa_1\kappa_2\kappa_3\neq 0$} By Proposition~\ref{crit}, $d_i\neq 0$ and thus $\kappa_i\neq \kappa_j$. We thus have the relations~(\ref{relations}), which we need to solve for~$\{u_i,v_i\}$ in~$\mathbf{Z}$. Two cases arise:

\subsubsection{When $u_iv_i\neq 6$ for every~$i$} If we impose this condition, all the integer solutions (in terms of~$u_i,v_i$) to equations~(\ref{relations}),  may be constructed. The first one, as a Diophantine equation in~$d_1$, $d_2$ and~$d_3$, has a finite number of solutions and all of them may be algorithmically  found (if~$d_1$ denotes the biggest~$d_i$, $0<d_1\leq 18$, and thus there are but finitely many choices for~$d_1$; repeating the argument, there are finitely many choices for the remaining~$d_i$). Each~$d_i$ in one of these solutions may be factored into potential~$u_i$ and~$v_i$ and, in turn, these may be evaluated in the remaining equations, in order to construct \emph{all} the solutions of the equations~(\ref{relations}) with~$u_i,v_i\in\mathbf{Z}$. In this case ($d_1\neq 6$ for every~$i$), we only have the solutions
$\{(2,5),(2,5),(-3,10)\}$ and $\{(2,5),(3,4),(-5,12)\}$. In these, $u_i+v_i=7$ for every~$i$ and hence~$\delta=0$. There is no longer one single solution for the~$\{\kappa_i\}$, but two of them. In the first case, these different solutions give a single solution for~$\{a_i\}$ and correspond to Chazy~IX. In the second one we still have two solutions for~$\{a_i\}$: they correspond to the two equations Chazy~X obtained by choosing a determination of~$\sqrt{3}$.

\subsubsection{When $u_1v_1= 6$} We necessarily have $(u_1,v_1)\in\{(2,3), (-2,-3), (1,6), (-1,-6)\}$. If~$(u_1,v_1)=(-2,-3)$ then the equations~(\ref{relations}) can be solved in terms of~$u_2=k$ and we get the family of solutions $\{(-2,-3), (1,k), (1,-k)\}$, which leads to Chazy XII if $k\notin\{1,6\}$. If~$k=1$ then~$L-\frac{5}{2}V$ vanishes at the point~$(1,\frac{2}{5},\frac{8}{25})$ where its linear part is
$$\left(\begin{array}{rrr} 1 & -\frac{5}{2} & 0 \\ 0 & 2 & -\frac{5}{2} \\  0 & \frac{6}{5} & -2 \end{array} \right),$$
but which is not diagonalizable (and thus~$V$ is not semicomplete). If $k=6$ there is no differential equation of the form (\ref{field:chazy}) that realizes it since~$k+1=7$ but~$1-k\neq 7$.

If~$(u_1,v_1)=(2,3)$ then the above equations can be solved in terms of~$u_2=k$ and we get the family of solutions $\{(2,3), (6,k), (6,-k)\}$ which leads to Chazy XI if~$k \neq 1$ (if~$k=1$ then $6+k=7$ but~$6-k\neq 7$ and thus there is no equation having this data).

If~$(u_1,v_1)=(1,6)$ then~$t_1=7$ and hence~$t_i=7$. From the first of the relations~(\ref{relations}), $d_2+d_3=0$ and we can assume, without loss of generality, that~$d_2>0$. Hence, since~$u_2+v_2=7$ and~$u_2v_2>0$, $(u_2,v_2)\in\{(1,6), (2,5), (3,4)\}$. This determines~$u_3$ and~$v_3$: they are not integers.

Finally, if~$(u_1,v_1)=(-1,-6)$, the second equation from~(\ref{relations}) reads~$u_2^{-1}+v_2^{-1}+u_3^{-1}+v_3^{-1}=14/6$. Since this number is greater than two, at least one of the variables must be equal to~$1$ (there is at least one positive summand greater than~$1/2$). This determines the remaining variables: they are never integers.

\subsection{Case~2, $\kappa_1=0$, $\kappa_2\kappa_3\neq 0$} By Proposition~\ref{crit}, $\kappa_2\neq\kappa_3$. We may still define~$t_2$, $t_3$, $d_2$ and $d_3$ like in~(\ref{defeigenvals}) and we have the two relations
$$\frac{1}{d_2}+\frac{1}{d_3}=\frac{1}{6},\; \frac{t_2}{d_2}+\frac{t_3}{d_3}=\frac{7}{6}.$$

The integer solutions to both equations in terms of the~$(u_i,v_i)$ are $\{(1,3),(-2,3)\}$ (leading to Chazy~IV), $\{(1,4),(-3,4)\}$ (Chazy~V), $\{(1,5),(-5,6)\}$ (Chazy~VI), $\{(2,4),(4,6)\}$ (Chazy~VII), $\{(3,4),(3,4)\}$ (Chazy~VIII; in this case, $\delta=0$ and there are two solutions at the level of the~$\kappa_i$ that become a single one at the level of the~$a_i$) and $\{(3,3),(3,6)\}$, leading to the equation~$\phi'''=3\phi^2\phi'+ (\phi')^2+ \phi\phi''$. For the latter, for the corresponding vector field~$V$ of the form~(\ref{field:chazy}), $L-V$ vanishes
at the point $(1,1,2)$. Its linear part at this point is
 $$\left(\begin{array}{rrr} 1 & -1 & 0 \\ 0 & 2 & -1 \\ -8 & -5 & 2 \end{array}\right).$$
The matrix has eigenvalues $-1$, $3$, $3$ but is not diagonalizable and the vector field cannot be semicomplete.

\subsection{Case~3, $\kappa_2=0, \kappa_3=0$, $\kappa_1\neq0$} We may still define $t_1$ and~$d_1$. We have~$d_1=6$ and thus~$t_1\in\{-7,-5,5,7\}$.
In the cases where~$(u_1,v_1)$ equals $(1,6)$, $(2,3)$, $(-2,-3)$ we obtain, respectively, Chazy~XI (in the case~$k=1)$, Chazy~II and Chazy~XII (in the case~$k=\infty)$. In the case $(-1,-6)$, the resulting equation is
$\phi'''=11(\phi')^2-7\phi\phi''$. For the corresponding vector field~$V$ of the form~(\ref{field:chazy}), $L-V$ vanishes at the point~$p=(1,-\frac{1}{2},\frac{1}{2})$, where its linear part is
$$\left(\begin{array}{rrr} 1 & 2 & 0 \\ 0 & 2 & 2 \\ -7 & -22 & -10 \end{array}\right).$$
This matrix has eigenvalues~$-1,-1,-6$ and is not diagonalizable: the vector field cannot be semicomplete.

\subsection{Case~4, $\kappa_i=0$ for every~$i$} In this case the equation is $\phi'''=\phi\phi''-2(\phi')^2$. The corresponding vector field vanishes along the curve~$t\mapsto(t,0,0)$ and, off this curve, is always linearly independent with~$L$. Our criterion does not yield any obstruction for semicompleteness. We refer the reader to~\cite{chazy-limitation} for the proof of the fact that not every solution of this equation is single-valued.

\end{document}